\documentclass[a4paper,12pt]{amsart}

\usepackage{amsmath}
\usepackage{amssymb}
\usepackage{mathrsfs}
\usepackage{ifthen}
\usepackage{graphicx}
\usepackage[T1]{fontenc} %skandit

\def\switchlinenumbers{\@ifstar
    {\let\makeLineNumberOdd\makeLineNumberRight
     \let\makeLineNumberEven\makeLineNumberLeft}%
    {\let\makeLineNumberOdd\makeLineNumberLeft
     \let\makeLineNumberEven\makeLineNumberRight}%
    }

\def\setmakelinenumbers#1{\@ifstar
  {\let\makeLineNumberRunning#1%
   \let\makeLineNumberOdd#1%
   \let\makeLineNumberEven#1}%
  {\ifx\c@linenumber\c@runninglinenumber
      \let\makeLineNumberRunning#1%
   \else
      \let\makeLineNumberOdd#1%
      \let\makeLineNumberEven#1%
   \fi}%
  }

\nonstopmode \numberwithin{equation}{section}
\newtheorem*{theorem*}{Theorem}

\newtheorem{thm}{Theorem}[section]
\newtheorem{cor}[equation]{Corollary}
\newtheorem{lem}[equation]{Lemma}
\newtheorem{prop}[equation]{Proposition}

\theoremstyle{definition}
\newtheorem{defn}{Definition}[section]

\newtheorem{prob}[equation]{Problem}

\newcounter{minutes}
\divide\time by 60
\newcounter{hours}
\multiply\time by 60
\addtocounter{minutes}{-\time}
\newcounter {own}
\def\theown {\thesection       .\arabic{own}}

\newenvironment{pf}[1][]{%
 \vskip 3mm
 \noindent
 \ifthenelse{\equal{#1}{}}%
  {{\slshape Proof. }}%
  {{\slshape #1.} }%
 }%
{\qed\bigskip}

\newcounter{alphabet}

\def\be{\begin{equation}}
\def\ee{\end{equation}}

\newcommand{\bee}{\begin{enumerate}}
\newcommand{\eee}{\end{enumerate}}

\newcommand{\blem}{\begin{lem}}
\newcommand{\elem}{\end{lem}}
\newcommand{\bthm}{\begin{thm}}
\newcommand{\ethm}{\end{thm}}
\newcommand{\bcor}{\begin{cor}}
\newcommand{\ecor}{\end{cor}}
\newcommand{\beg}{\begin{examp}}
\newcommand{\eeg}{\end{examp}}
\newcommand{\begs}{\begin{examples}}
\newcommand{\eegs}{\end{examples}}
\newcommand{\bdefe}{\begin{defin}}
\newcommand{\edefe}{\end{defin}}
\newcommand{\bprob}{\begin{prob}}
\newcommand{\eprob}{\end{prob}}
\newcommand{\bei}{\begin{itemize}}
\newcommand{\eei}{\end{itemize}}

\newcommand{\norm}[1]{\left\lVert#1\right\rVert}

\newcommand{\innpdct}[1]{\left\langle#1\right\rangle}

\begin{document}

\title{On multidimensional Bohr radius of finite dimensional Banach spaces}

\author{Vasudevarao Allu}
\address{Vasudevarao Allu,
Department of Mathematics,
School of Basic Sciences,
Indian Institute of Technology Bhubaneswar,
Bhubaneswar-752050, Odisha, India.}
\email{avrao@iitbbs.ac.in}

\author{Subhadip Pal}
\address{Subhadip Pal,
	Department of Mathematics,
	School of Basic Sciences,
	Indian Institute of Technology Bhubaneswar,
	Bhubaneswar-752050, Odisha, India.}
\email{subhadippal33@gmail.com}

\subjclass[{AMS} Subject Classification:]{Primary 32A05, 32A10; Secondary 46B07, 46G25}

\def\thefootnote{}
\footnotetext{ {\tiny File:~\jobname.tex,
printed: \number\year-\number\month-\number\day,
          \thehours.\ifnum\theminutes<10{0}\fi\theminutes }
} \makeatletter\def\thefootnote{\@arabic\c@footnote}\makeatother

\begin{abstract}
In this paper, we improve the lower estimate of multidimensional Bohr radius for unit ball of $\ell^n_q$-spaces ($1\leq q\leq \infty$) for bounded holomorphic functions with values in finite dimensional complex Banach spaces. The new estimate provides the improved lower bound for the Bohr radius which was previously given by Defant, Maestre, and Schwarting [Adv. Math. 231 (2012), 2837--2857].
\end{abstract}

\maketitle
\pagestyle{myheadings}
\markboth{Vasudevarao Allu and Subhadip Pal}{On multidimensional Bohr radius of finite dimensional Banach spaces}

\section{Introduction}
A complete Reinhardt domain $\Omega \subseteq \mathbb{C}^n$ is a domain in $\mathbb{C}^n$ such that if $z=(z_{1},\ldots,z_{n}) \in \Omega$, then $(\lambda e^{i\theta_{1}}z_{1},\ldots,\lambda e^{i\theta_{n}}z_{n}) \in \Omega$ for all $\lambda \in \overline{\mathbb{D}}$ and all $\theta_{i} \in \mathbb{R}$, $i=1, \ldots,n$. For a complete Reinhardt domain $\Omega\subset \mathbb{C}^n$, the Bohr radius $K(\Omega)$ is the supremum of all $r\geq 0$ such that for each holomorphic function $f(z)=\sum_{\alpha \in \mathbb{N}^{n}_{0}}c_{\alpha}z^{\alpha}$ on $\Omega$, we have 
\begin{equation*}
	\sup_{z\in r\Omega}\sum_{\alpha \in \mathbb{N}^{n}_{0}}|c_{\alpha}z^{\alpha}|\leq \sup_{z\in \Omega}\bigg|\sum_{\alpha \in \mathbb{N}^{n}_{0}}c_{\alpha}z^{\alpha}\bigg|.
\end{equation*}
Here, $\alpha$ denotes an $n$-tuple $(\alpha_1,\ldots,\alpha_n)\in \mathbb{N}^n_{0}$ ($\mathbb{N}_{0}=\mathbb{N}\cup\{0\}$) of non-negative integers, $z$ stands for an $n$-tuple $(z_1,\ldots,z_n)$ of complex numbers, and $z^{\alpha}$ denotes the product $z^{\alpha_1}_{1}\cdots z^{\alpha_n}_{n}$. By a remarkable result of Harald Bohr \cite{Bohr-1914}, also known as famous Bohr's power series theorem shows that $K(\mathbb{D})=1/3$. In this paper, we keep our interest on the complete Reinhardt domains in $n$-dimensions are
\begin{equation*}
	B_{\ell^n _q}:=\left\{z\in \mathbb{C}^n: \norm{z}_{q}<1\right\}, \quad 1\leq q\leq \infty,
\end{equation*}
where we write $\ell^n_q$ for the Banach space defined by $\mathbb{C}^n$ together with the $q$-norm $\norm{z}_q:=\left(\sum_{i=1}^{n}|z_{i}|^q\right)^{1/q}$, for $1\leq q<\infty$ and $\norm{z}_{\infty}:=\sup\{|z_{i}|:1\leq i\leq n\}$ for $q=\infty$. In 1997, Boas and Khavinson \cite{boas-1997} studied the $n$ dimensional Bohr radius $K^n(B_{\ell^n_{\infty}})$ for the polydisk $B_{\ell^n_{\infty}}$ in $\mathbb{C}^n$ and obtained the following interesting estimate:
\begin{equation}\label{P7-e-011}
	\frac{1}{3\sqrt{n}}\leq K^n(B_{\ell^n_q})\leq \frac{2\sqrt{\log n}}{\sqrt{n}}.
\end{equation}
However, the concept of Bohr radius in several variables first appeared implicitly in a paper of Dineen and Timoney \cite{Dineen-Timoney-1989} in which the authors have explained the idea of extending the Bohr's result in higher dimension. Later, the problem becomes very active in regard to find the correct asymptotic estimates of $n$ dimensional Bohr radius $K^n(B_{\ell^n_{q}})$ for $1\leq q\leq \infty$. Investigations by the galaxy of analysts like Aizenberg, Boas, Defant, Dineen, Frerick, Khavinson and Timoney from \cite{aizn-2000a}, \cite{boas-2000}, \cite{boas-1997}, \cite{defant-2006}, \cite{Dineen-Timoney-1989} established that there is a constant $C\geq 1$ such that for all $1\leq q\leq \infty$, and all $n\in \mathbb{N}$
\begin{equation}\label{P7-e-010}
	\frac{1}{C}\left(\frac{\log n}{n\log \log n}\right)^{1-\frac{1}{\min\{q,2\}}}\leq K^n(B_{\ell^n_q})\leq C\left(\frac{\log n}{n}\right)^{1-\frac{1}{\min\{q,2\}}}.
\end{equation}
By use of surprising fact of hypercontarctivity of Bohnenblust-Hille inequality for homogeneous polynomials, a groundbreaking progress on this problem, especially for the case $q=\infty$ has made by Dafant and his coauthors in \cite[Theorem 2]{defant-2011} and they have shown that
\begin{equation*}
	K^n(B_{\ell^n_{\infty}})=b_n\sqrt{\frac{\log n}{n}} \quad \mbox{with} \quad \frac{1}{\sqrt{2}} + o(1)\leq b_n\leq 2.
\end{equation*}
 Finally, Bayart {\it et. al.} \cite{bayart-advance-2014} have established the exact asymptotic behaviour of $K^n(B_{\ell^n_{\infty}})$ as
 \begin{equation*}
 	K^n(B_{\ell^n_{\infty}})\sim_{+\infty} \sqrt{\frac{\log n}{n}}.
 \end{equation*}
Numerous significant contributions have been made to the study of the multidimensional Bohr radius. Notably, the works of Aizenberg and his collaborators \cite{aizn-2000a,aizn-2000b,aizenberg-2001,aizn-2005} have explored this topic in the context of various domains in $\mathbb{C}^n$. Important developments have also been made by Balasubramanian, Calado, and Queffélec \cite{bala-studia-2006} in the setting of Dirichlet series; by Defant and his coauthors \cite{defant-2003,defant-2011-angew,defant-2012,defant-2018} in the study of vector-valued holomorphic functions; by Paulsen and Singh \cite{paulsen-2002} in the framework of Banach algebras; and by Popescu \cite{popescu-2019} in the context of free holomorphic functions.\\

The study of Bohr radius becomes very interesting when one considers the space of holomorphic functions in Reinhardt domains in $\mathbb{C}^n$ with values in Banach spaces (vector-valued) instead of scalar-valued functions. Initially, for the one dimensional case, the Bohr radius $K(\mathbb{D},X)$ for the space of holomorphic functions from unit disk $\mathbb{D}$ in complex plane to any arbitrary complex Banach space $X$ was studied by Blasco \cite{Blasco-OTAA-2010} and followed by $p$-Bohr radius $K_p(\mathbb{D},X)$ for powered Bohr's inequality for any $p\in [1,\infty]$ can be found in \cite{Blasco-Collect-2017}. It has been noticed by Blasco \cite{Blasco-OTAA-2010} that for $X=\ell^2_q$ the Bohr radius $K(\mathbb{D},\ell^2_q)=0$ for all $1\leq q\leq \infty$. This fact motivates Defant {\it et. al.} \cite{defant-2012} to modify the definition of Bohr radius in order to obtain non-trivial Bohr radius. Indeed, Defant {\it et. al.} \cite{defant-2012} have introduced the Bohr radius in more general form.
\begin{defn}\cite{defant-2012}
	Let $T:X\rightarrow Y$ be a bounded linear operator between complex Banach spaces, $n\in \mathbb{N}$, and $\norm{T}\leq \lambda$. The $\lambda$-Bohr radius of $T$, denoted by $K^n(B_{\ell^n_{\infty}},T,\lambda)$ is the supremum of all $r\geq 0$ such that for all holomorphic functions $f(z)=\sum_{\alpha \in \mathbb{N}^{n}_{0}}x_{\alpha}z^{\alpha}$ on $B_{\ell^n_{\infty}}$ we have 
	\begin{equation*}
		\sup_{z\in rB_{\ell^n_{\infty}}} \sum_{\alpha \in \mathbb{N}^{n}_{0}}\norm{T(x_{\alpha})z^{\alpha}}\leq \lambda\sup_{z\in B_{\ell^n_{\infty}}} \norm{\sum_{\alpha \in \mathbb{N}^{n}_{0}}x_{\alpha}z^{\alpha}}.
	\end{equation*}
\end{defn}
\noindent
 If $T$ is the identity operator on $X$ we denote the $\lambda$-Bohr radius as $K^n(B_{\ell^n_{\infty}},X,\lambda)$. For $\lambda=1$ and $X=\mathbb{C}$, we use the notation $K^n(B_{\ell^n_{\infty}},X)$ and $K^n(B_{\ell^n_{\infty}})$ respectively.\\ 
 
 In \cite{defant-2012}, the authors have studied the Bohr radius $K^n(B_{\ell^n_{\infty}},X,\lambda)$ for both finite and infinite dimensional complex Banach spaces. In particular for finite dimensional complex Banach spaces, Dafant and his co-authors have proved the following:
 \begin{thm}
 	Let $X$ be a finite dimensional complex Banach space and $\lambda>1$. Then there are constants $C$, $D>0$ depending on $\lambda$ and $X$ such that for each $n\in \mathbb{N}$
 	\begin{equation}\label{P7-e-008}
 		C\sqrt{\frac{\log n}{n}}\leq K^n(B_{\ell^n_{\infty}},X,\lambda)\leq D\sqrt{\frac{\log n}{n}}.
 	\end{equation}
 \end{thm}
\noindent
 Recently, Kumar and Manna \cite{kumar-2023-arxiv}, have estimated the Bohr radius $K^n(B_{\ell^n_{q}},X,\lambda)$ for the remaining unit ball $B_{\ell^n_q}$, $1\leq q<\infty$ and proved that for each $1\leq q<\infty$ and $\lambda>1$,
 \begin{equation}\label{P7-e-009}
 	C \left(\frac{\log n}{n}\right)^{1-\frac{1}{\min\{q,2\}}}\leq K^n(B_{\ell^n_{q}},X,\lambda) \leq D\left(\frac{\log n}{n}\right)^{1-\frac{1}{\min\{q,2\}}},
 \end{equation}
 where the constants $C$ and $D$ depends on $\lambda$ and $X$. The main aim of this paper is to show an improved lower bound of \eqref{P7-e-008} and \eqref{P7-e-009}. In 2021, Bernal-Gonz\'{a}lez {\it et. al.} implicitly \cite[Theorem 3.3]{Bernal-2021} improved the lower bound of $K^n(B_{\ell^n_{\infty}})$ given in \eqref{P7-e-011}. Indeed, the authors \cite{Bernal-2021} have proved that the Bohr radius $K^n(B_{\ell^n_{\infty}})$ is not less than the solution $r>0$ of the equation
 \begin{equation*}
 	r+\sum_{m=2}^{\infty}S(m,n)r^m=\frac{1}{2},
 \end{equation*}
where $S(m,n)$ is the Sidon constant for each $m,n\in \mathbb{N}$ and it is defined by the infimum of all $\eta>0$ such that 
\begin{equation}\label{P7-e-013}
	\sum_{\alpha \in \Lambda(m,n)}|c_{\alpha}|\leq \eta \sup_{z\in B_{\ell^n_{\infty}}}\bigg|\sum_{\alpha \in \Lambda(m,n)}c_{\alpha}z^{\alpha}\bigg|
\end{equation} 
for all $m$-homogeneous polynomials $\sum_{\alpha \in \Lambda(m,n)}c_{\alpha}z^{\alpha}\in \mathcal{P}(^m\ell^n_q,\mathbb{C})$. Recently, Kumar and Manna \cite{kumar-2023-JMAA} demonstrated that the methods introduced in \cite{Bernal-2021} can be effectively applied to improve the lower bound of $K^n(B_{\ell^n_q})$ for $1 \leq q < \infty$, as presented in \eqref{P7-e-010}. Our aim is to implement the idea for vector-valued holomorphic functions and establish a new lower bound for $K^n(B_{\ell^n_q},X,\lambda)$ for any finite dimensional complex Banach spaces $X$ and all $1\leq q\leq \infty$, $\lambda>1$.
In Section \ref{P7-sec-02}, we shall provide basic definitions, notations, and some preliminary statements. In Section \ref{P7-sec-03}, we shall discuss the connection between Bohr radius and unconditional basis constant for space of vector-valued homogeneous polynomials, and the main results of this paper is contained in the final section i.e., in Section \ref{P7-sec-04}.
\section{Preliminary results and Notations}\label{P7-sec-02}
Following the notation from \cite{Bernal-2021}, for every $n,m\in \mathbb{N}$, we write $N_m(n)$ as the combinatorial number 
\begin{equation*}
	N_m(n):=\binom{n+m-1}{m}=\frac{(n+m-1)!}{(n-1)!m!}.
\end{equation*}
Using the Cauchy-Hadamard formula one can easily find that the radius of convergence of the series $\sum_{m=2}^{\infty}\sqrt{N_m(n)}x^m$ is $1$. Therefore, the function 
\begin{equation}\label{P7-e-012}
	H_n(x):=x+\sum_{m=2}^{\infty}\sqrt{N_m(n)}x^m
\end{equation}
is well-defined and differentiable in $(-1,1)$, and so in $[0,1)$. However, we have the following properties of $H_n$ discussed in \cite[Lemma 2.1]{Bernal-2021}.
\begin{lem}\cite{Bernal-2021}\label{P7-lem-04}
	For every $n \in \mathbb{N}$, let $H_n: [0,1) \rightarrow \mathbb{R}$ be the restriction to $[0,1)$ of the function defined by equation \eqref{P7-e-012}. Then we have
	\begin{enumerate}
		\item[(a)] $H'_{n}(x)>0$ for all $x\in (0,1)$.
		\item[(b)] $H_n$ is strictly increasing.
		\item[(c)] $H_n$ is injective.
		\item[(d)] $\lim_{x\rightarrow 1^{-}}H_n(x)=+\infty$.
		\item[(e)] For every $M\in (0,+\infty)$ there exists a unique $S\in (0,1)$ such that $H_n(S)=M$.
	\end{enumerate}
\end{lem}
\noindent
In view of Lemma \ref{P7-lem-04}(e), it is evident that the following definition makes sense.

\begin{defn}
	For every $n\in \mathbb{N}$ and $1\leq \lambda<\infty$, we denote by $\beta_n \in (0,1)$ the unique solution of the following equation:
	\begin{equation}\label{P7-e-001-beta}
		x +\sum_{m=2}^{\infty}\sqrt{N_m(n)}x^m=\frac{\lambda}{2}.
	\end{equation}
\end{defn}

Throughout our discussion all Banach spaces $X$ are assumed to be complex and finite dimensional. The topological duals of $X$ are denoted by $X^*$ and we write their unit balls as $B_{X}$. For each $n\in \mathbb{N}$ and $\alpha=(\alpha_1,\ldots,\alpha_n)\in \mathbb{N}^n_{0}$, we use the standard multi-index notation $|\alpha|=\alpha_1+\cdots+\alpha_n$. For every $m,n\in \mathbb{N}$, the symbol $\Lambda(m,n)$ stands for the set of all $\alpha=(\alpha_1,\ldots,\alpha_n)\in \mathbb{N}^n_{0}$ such that $|\alpha|=m$.
Let $\mathcal{P}(^m\ell^n_q,X)$ denote the finite dimensional Banach space of all vector-valued $m$-homogeneous polynomials $Q(z)=\sum_{\alpha\in \Lambda(m,n)}x_{\alpha}z^{\alpha},$ where $x_{\alpha}\in X$ and $z\in \mathbb{C}^n$. For every $Q(z)=\sum_{\alpha\in \Lambda(m,n)}x_{\alpha}z^{\alpha}\in \mathcal{P}(^m\ell^n_q,X)$, the norms $\norm{\cdot}_{\infty}$ and $\norm{\cdot}_{1}$ are defined as follows
\begin{equation*}
	\norm{Q}_{\infty}:=\sup_{z\in \partial B_{\ell^n_q}}\norm{Q(z)} \quad \mbox{and}\quad \norm{Q}_{1}:=\sup_{z\in \partial B_{\ell^n_q}}\sum_{\alpha\in \Lambda(m,n)}\norm{x_{\alpha}z^{\alpha}},
\end{equation*}
where $\partial B_{\ell^n_q}=\{z\in B_{\ell^n_q}:\norm{z}_q=1\}$. Clearly, $\norm{Q}_{\infty}\leq \norm{Q}_{1}$.\\

For every $z\in \ell^n_q$, let us denote $\Omega_z:= \{(z_1e^{2\pi i\theta_{1}}, \ldots, z_{n}e^{2\pi i\theta_n}): \theta_{i}\in [0,1], 1\leq i\leq n\}$. Then 
\begin{align*}
	\norm{Q(z)}_2=\left(\sum_{\alpha\in \Lambda(m,n)}\norm{x_{\alpha}z^{\alpha}}^2\right)^{1/2} 
	&= \left(\int_{[0,1]^n} \norm{Q(z_1 e^{2\pi i\theta_{1}}, \ldots, z_n e^{2\pi i\theta_{n}})}^2 \,d\theta_{1}\cdots d\theta_{n}\right)^{1/2}\\ &
	\leq \sup_{v\in \Omega_{z}}\norm{Q(v)}_{\mathcal{P}(^m\ell^n_q,X)}=\sup_{\substack{x^*\in X^* \\ \norm{x^*}\leq 1}} \sup_{v\in \Omega_{z}}|x^*(Q(v))|.
\end{align*}
Moreover, it is not so difficult to see that $\norm{Q(z)}_2=\sup_{v\in \Omega_{z}}\norm{(Q(v))}$ if, and only if, $|x^*(Q(v))|$ is constant for every $x^*\in X^*$ with $\norm{x^*}\leq 1$.\\

For each $j\in \mathbb{N}$, let us denote by $e_{j}$ the $n$-tuple $(0,0,\ldots,0,1,0,\ldots,0)$ with $1$ at the $j$th place. Further, $me_{j}$ will denote the $n$-tuple with $m$ at the $j$th place, and zero elsewhere.

\begin{lem}\label{P7-lem-002}
	Assume that $X$ be any finite dimensional complex Banach space and $n\geq 2$, $1\leq q\leq \infty$. Let $Q\in \mathcal{P}(^m\ell^n_q,X)$ with 
	\begin{equation*}
		Q(z)=\sum_{\alpha\in \Lambda(m,n)}x_{\alpha}z^{\alpha}
	\end{equation*}
	satisfying that there exist $j,j'$ and $x^*\in B_{X^*}$ such that $x^*(z_j x_{me_{j}})\neq 0 \neq x^*(z_{j'}x_{me_{j'}})$ for $j\neq j'$. Then 
	\begin{equation*}
		\norm{Q(z)}_{2}< \sup_{v\in \Omega_{z}}\norm{Q(v)}_{\mathcal{P}(^m\ell^n_q,X)}.
	\end{equation*}
	
\end{lem}
\begin{pf}
	For $Q\in \mathcal{P}(^m\ell^n_q,X)$ and $x^* \in B_{X^*}$, we define the trigonometric polynomial $T_{Q}$ of degree $2m$ and $n$ variables as
	\begin{equation*}
		T_{Q}(\widetilde{\theta}):= |x^*(Q(z_1 e^{i\theta_1},\ldots,z_ne^{i\theta_n}))|^2=x^*(Q(z_1 e^{i\theta_1},\ldots,z_ne^{i\theta_n}))\overline{x^*(Q(z_1 e^{i\theta_1},\ldots,z_ne^{i\theta_n}))}
	\end{equation*}
	for all $\widetilde{\theta}=(\theta_{1},\ldots,\theta_{n})\in \mathbb{R}^n$ and $z=(z_1, \ldots,z_n)\in \partial B_{\ell^n_q}$. If we use the notation $\innpdct{\alpha, \widetilde{\theta}}=\sum_{k=1}^{n}\alpha_k \theta_k$, then we have that
	\begin{align*}
		\frac{\partial T_{Q}}{\partial \theta_{j}}&=i\sum_{\alpha\in \Lambda(m,n)}\alpha_j x^*(x_{\alpha})z^{\alpha}e^{i\innpdct{\alpha, \widetilde{\theta}}} \sum_{\alpha\in \Lambda(m,n)}\overline{x^*(x_{\alpha})}\overline{z}^{\alpha}e^{-i\innpdct{\alpha, \widetilde{\theta}}}\\ & -i\sum_{\alpha\in \Lambda(m,n)} x^*(x_{\alpha})z^{\alpha}e^{i\innpdct{\alpha, \widetilde{\theta}}} \sum_{\alpha\in \Lambda(m,n)}\alpha_j\overline{x^*(x_{\alpha})}\overline{z}^{\alpha}e^{-i\innpdct{\alpha, \widetilde{\theta}}}\\ &
		=i\sum_{\substack{\alpha, \alpha'\in \Lambda(m,n)\\\alpha\neq \alpha'}} \alpha_j x^*(x_{\alpha})z^{\alpha}\overline{x^*(x_{\alpha'})}\overline{z}^{\alpha'}e^{i\innpdct{\alpha-\alpha', \widetilde{\theta}}}-i\sum_{\substack{\alpha, \alpha'\in \Lambda(m,n)\\\alpha\neq \alpha'}} \alpha_j x^*(x_{\alpha'})z^{\alpha'}\overline{x^*(x_{\alpha})}\overline{z}^{\alpha}e^{i\innpdct{\alpha'-\alpha, \widetilde{\theta}}}\\ &
		=imx^*(x_{me_{j}})\overline{x^*(x_{me_{j'}})}z^m_{j} z^m_{j'}e^{im(\theta_{j}-\theta_{j'})}+\cdots-imx^*(x_{me_{j'}})\overline{x^*(x_{me_{j}})}z^m_{j'} z^m_{j}e^{im(\theta_{j'}-\theta_{j})}-\cdots
	\end{align*}
	where $\alpha-\alpha'=(\alpha_1-\alpha_1',\dots,\alpha_n-\alpha_n')$.
	Since $x^*(x_{me_{j}})x^*(x_{me_{j'}})\neq 0$, it follows that 
	\begin{equation*}
		\frac{\partial T_{Q}}{\partial \theta_{j}}\neq 0,
	\end{equation*}
	which shows that $T_{Q}$ is non-constant. Further, $T_{Q}$ is merely $|x^*(Q)|^2$. Therefore, we obtain $|x^*(Q)|$ is non-constant. This competes the proof.
\end{pf}

For every $Q\in \mathcal{P}(^m\ell^n_q,X)$, we define $\norm{Q}_2=\sup_{z\in \partial B_{\ell^n_q}}\norm{Q(z)}_2$. Then, using the fact card $\Lambda(m,n)=\mbox{card}\,\{\alpha \in \mathbb{N}^{n}_{0}: |\alpha|=m\}=\binom{n+m-1}{m}=N_m(n)$ we obtain
\begin{equation}\label{P7-e-003}
	\norm{Q}_1=\sum_{|\alpha|=m}\norm{x_{\alpha}z^{\alpha}}\leq \left(\sum_{|\alpha|=m}1\right)^{1/2} \left(\sum_{|\alpha|=m}\norm{x_{\alpha}z^{\alpha}}^2\right)^{1/2}=\sqrt{N_m(n)} \norm{Q}_2.
\end{equation}

\begin{prop}\label{P7-prop-01}
	Let $n\in \mathbb{N}$ with $n\geq 2$ and $X$ be any finite dimensional complex Banach space. Then for all $Q\in \mathcal{P}(^m\ell^n_q,X)$ with $\norm{Q}_{\infty}=1$ we have 
	\begin{equation*}
		\norm{Q}_1< \sqrt{N_m(n)}.
	\end{equation*}
\end{prop}
\begin{pf}
	We first note that $\norm{Q}_2\leq \norm{Q}_{\infty}$. Then, in view of \eqref{P7-e-003}, we have that 
	\begin{equation}\label{P7-e-004}
		\norm{Q}_1\leq \sqrt{N_m(n)} \norm{Q}_2\leq \sqrt{N_m(n)}.
	\end{equation}
	We prove the required assertion by the method of contradiction. Let us suppose there exists $Q(z)=\sum_{\alpha\in \Lambda(m,n)}v_{\alpha}z^{\alpha}\in \mathcal{P}(^m\ell^n_q,X)$ with $\norm{Q}_{\infty}=1$ such that $\norm{Q}_{1}=\sqrt{N_m(n)}$. Then, by the virtue of the inequality \eqref{P7-e-004}, it follows that 
	\begin{equation*}
		\norm{Q}_{1}=\sqrt{N_m(n)}\norm{Q}_{2}=\sqrt{N_m(n)}.
	\end{equation*}
	Hence, we obtain $\norm{Q}_{2}=\norm{Q}_{\infty}=1$, which shows that
	\begin{equation*}
		\norm{x_{\alpha}\widetilde{z}^{\alpha}}=\frac{\left(\sum_{\alpha\in \Lambda(m,n)}\norm{x_{\alpha}\widetilde{z}^{\alpha}}^2\right)^{1/2}}{\sqrt{N_m(n)}}=\frac{1}{\sqrt{N_m(n)}}.
	\end{equation*}
	Finally, in view of Lemma \ref{P7-lem-002}, we conclude that 
	\begin{equation*}
		1=\norm{Q(\widetilde{z})}_{2}<\sup_{v\in \Omega_{z}}\norm{Q(v)}\leq 1,
	\end{equation*}
	which is a contradiction. This completes the proof.
\end{pf}

\begin{defn}\label{defn-sidon}
	For each pair $m,n\in \mathbb{N}$, we define the constant $\widetilde{S}(m,n)$ as 
	\begin{align*}
		\widetilde{S}(m,n)&=\sup\left\{\norm{Q}_1: Q(z)=\sum_{\alpha\in \Lambda(m,n)}x_{\alpha}z^{\alpha}\in \mathcal{P}(^m\ell^n_q,X)\,\, \mbox{and}\,\, \norm{Q}_{\infty}\leq 1  \,\,\mbox{for all}\,\, z\in B_{\ell^n_q}\right\}\\[2mm]&
		=\inf \left\{\eta>0: \norm{Q}_1\leq \eta \norm{Q}_{\infty} \,\,\mbox{for all}\,\, P\in \mathcal{P}(^m\ell^n_q,X)\right\}.
	\end{align*}
\end{defn}
\noindent
It is worth noting that for $X=\mathbb{C}$ and $q=\infty$, the constant $\widetilde{S}(m,n)$ can be viewed as the usual Sidon constant $S(m,n)$ defined in \eqref{P7-e-013}. Evidently, $\widetilde{S}(m,n)\geq 1$ for every $m,n\in \mathbb{N}$. Further, observe that the constant $\widetilde{S}(m,n)$ coincides with the norm of the identity operator
\begin{equation*}
	I:\left(\mathcal{P}(^m\ell^n_q,X),\norm{\cdot}_{\infty}\right) \rightarrow \left(\mathcal{P}(^m\ell^n_q,X), \norm{\cdot}_{1}\right).
\end{equation*}
Since the space $\mathcal{P}(^m\ell^n_q,X)$ is finite dimensional, for every $m\in \mathbb{N}$ we can find $Q_m\in \mathcal{P}(^m\ell^n_q,X)$ satisfying 
$\norm{Q}_{\infty}=1$ and $\norm{Q}_{1}=\widetilde{S}(m,n)$. Consequently, in view of Proposition \ref{P7-prop-01}, we obtain the following corollary.
\begin{cor}\label{P7-cor-01}
	For each pair $m,n\in \mathbb{N}$ with $n\geq 2$, we have 
	\begin{equation*}
		0<\widetilde{S}(m,n)<\sqrt{N_{m}(n)}.
	\end{equation*}
\end{cor}

\section{Bohr radius and unconditional basis constant}\label{P7-sec-03}
A Schauder basis $\{x_n\}_{n\geq 1}$ of a Banach space $E$ is said to be unconditional if there is a constant $C\geq 1$ such that 
\begin{equation*}
	\norm{\sum_{j=1}^{n}\epsilon_j \sigma_j x_{j}}\leq C \norm{\sum_{j=1}^{n}\sigma_j x_{j}}	
\end{equation*}
for all $n\in \mathbb{N}$, and all $\sigma_1, \ldots, \sigma_n \in \mathbb{C}$, and all $\epsilon_1, \ldots,\epsilon_n \in \mathbb{C}$ with $|\epsilon_k|\leq 1$, $1\leq k\leq n$. In this case, the best constant $C$ is denoted by $\chi(\{x_n\})$ and called the unconditional basis constant of $\{x_n\}$. Moreover, the unconditional basis constant $\chi(E)\in [1, \infty]$ of the Banach space $E$ is defined by the infimum over all constants $\chi(\{x_n\})$. We say that a basis $\{x_n\}$ of a Banach space is $1$-unconditional if $\chi(\{x_n\})=1$. In this paper, we mainly consider the $n$-dimensional Banach spaces $X=(\mathbb{C}^n, \norm{\cdot}_q)$, $1\leq q\leq \infty$ such that the canonical basis vectors $\{e_n\}_{n\geq 1}$ form a $1$-unconditional basis, or another way it means that the unit ball $B_{X}$ is a Reinhardt domain.\\

By $\chi_{\text{mon}}(\mathcal{P}(^m\ell^n_q,X))$, we denote the unconditional basis constant of the monomials $z^{\alpha}, \, |\alpha|=m$ in the Banach space $\mathcal{P}(^m\ell^n_q,X)$. Equivalently, we can write
\begin{equation*}
	\chi_{\text{mon}}(\mathcal{P}(^m\ell^n_q,X)):=\sup\left\{\norm{\sum_{\alpha\in \Lambda(m,n)}\norm{x_{\alpha}}z^{\alpha}}: \norm{\sum_{\alpha\in \Lambda(m,n)}x_{\alpha}z^{\alpha}}\leq 1\right\}.
\end{equation*}
Let $X$ be any finite dimensional Banach space. It is worth noting that the Bohr radius $K^n(B_{\ell^n_q},X,\lambda)$ of the open unit ball $B_{\ell^n_q}$ with respect to bounded  holomorphic functions with values in $X$ is the supremum of all $r\in [0,1]$ such that 

\begin{equation*}
	\sup_{z\in B_{\ell^n_q}}\sum_{\alpha \in \mathbb{N}^{n}_{0}}\norm{x_{\alpha}(rz)^{\alpha}}\leq \lambda \sup_{z\in B_{\ell^n_q}}\norm{\sum_{\alpha \in \mathbb{N}^{n}_{0}}x_{\alpha}z^{\alpha}}
\end{equation*}
for every power series expansion $\sum_{\alpha \in \mathbb{N}^{n}_{0}}x_{\alpha}z^{\alpha}$ convergent on $B_{\ell^n_q}$. Now, we have the following analogous Bohr radius for $m$-homogeneous polynomials.
\begin{defn}
	Let $1\leq q\leq \infty$ and $\lambda \geq 1$. For each $m\in \mathbb{N}$, we deine the Bohr radius, denoted by $K^n_m(B_{\ell^n_q},X,\lambda)$ and defined by the supremum of all $r\in [0,1]$ such that 
	\begin{equation*}
		\sup_{z\in rB_{\ell^n_q}}\sum_{\alpha \in \Lambda(m,n)}\norm{x_{\alpha}z^{\alpha}}\leq \lambda \sup_{z\in B_{\ell^n_q}} \norm{\sum_{\alpha\in \Lambda(m,n)}x_{\alpha}z^{\alpha}}
	\end{equation*}
for all $m$-homogeneous polynomials $\sum_{\alpha\in \Lambda(m,n)}x_{\alpha}z^{\alpha}$.
\end{defn}
Indeed, in view of triangle inequality and the fact that $\{|z_{i}|\}^n_{i=1}\in B_{\ell^n_q}$ for each $z\in B_{\ell^n_q}$, we have for each $m$ and each $m$-homogeneous polynomials $\sum_{\alpha\in \Lambda(m,n)}x_{\alpha}z^{\alpha}$ that 

\begin{align}
	\sup_{z\in B_{\ell^n_q}}\norm{\sum_{\alpha\in \Lambda(m,n)}\norm{x_{\alpha}}z^{\alpha}}&=\sup_{z\in B_{\ell^n_q}}\sum_{\alpha\in \Lambda(m,n)}\norm{x_{\alpha}z^{\alpha}}\\  \nonumber &\leq \frac{\lambda}{(K^n_m(B_{\ell^n_q},X,\lambda))^m}\sup_{z\in B_{\ell^n_q}} \norm{\sum_{\alpha\in \Lambda(m,n)}x_{\alpha}z^{\alpha}}.
\end{align}
It is not so hard to observe that $K^n(B_{\ell^n_q},X,\lambda) \leq K^n_m(B_{\ell^n_q},X,\lambda)$ for all $m$. Further, we note that $K^n_m(B_{\ell^n_q},X,\lambda)=\sqrt[m]{\lambda}K^n_m(B_{\ell^n_q},X,1)$.

\begin{lem}\label{P7-lem-03}
	Let $X$ be a finite dimesional complex Banach space. Then, for each $m\in \mathbb{N}$ and $\lambda \geq 1$ we have
	\begin{equation}\label{P7-e-015}
		K^n_m(B_{\ell^n_q},X,\lambda)=\frac{\sqrt[m]{\lambda}}{\sqrt[m]{	\chi_{\text{mon}}(\mathcal{P}(^m\ell^n_q,X))}}.
	\end{equation}
\end{lem} 

\begin{pf}
 We first prove the upper estimate for $K^n_m(B_{\ell^n_q},X,\lambda)$. Let $\sum_{\alpha\in \Lambda(m,n)}x_{\alpha}z^{\alpha}\in \mathcal{P}(^m\ell^n_q,X)$ be an $m$-homogeneous polynomial and $\epsilon_{\alpha}\in \mathbb{C}$ such that $|\epsilon_{\alpha}|\leq 1$. Then, we have 
 \begin{align*}
 	\norm{\sum_{\alpha\in \Lambda(m,n)}\epsilon_{\alpha}x_{\alpha}z^{\alpha}}_{\mathcal{P}(^m\ell^n_q,X)} &\leq \sup_{z\in B_{\ell^n_q}}\sum_{\alpha\in \Lambda(m,n)}\norm{x_{\alpha}z^{\alpha}}\\ &
 	\leq \frac{\lambda}{(K^n_m(B_{\ell^n_q},X,\lambda))^m}\sup_{z\in B_{\ell^n_q}} \norm{\sum_{\alpha\in \Lambda(m,n)}x_{\alpha}z^{\alpha}}.
 \end{align*}
Using the definition of unconditional basis constant of spaces of $m$-homogeneous polynomials, it follows that
\begin{equation*}
	\chi_{\text{mon}}(\mathcal{P}(^m\ell^n_q,X))\leq \frac{\lambda}{(K^n_m(B_{\ell^n_q},X,\lambda))^m},
\end{equation*}
which gives
\begin{equation}\label{P7-e-005}
	K^n_m(B_{\ell^n_q},X,\lambda)\leq \frac{\sqrt[m]{\lambda}}{\sqrt[m]{	\chi_{\text{mon}}(\mathcal{P}(^m\ell^n_q,X))}}.
\end{equation}
Conversely, given an $m$-homogeneous polynomial $\sum_{\alpha\in \Lambda(m,n)}v_{\alpha}z^{\alpha}$, where $v_{\alpha}\in X$ and $z\in \ell^n_q$, we obtain
\begin{align*}
	\sup_{z\in B_{\ell^n_q}}\sum_{\alpha\in \Lambda(m,n)}\norm{v_\alpha z^{\alpha}}&=\sup_{z\in B_{\ell^n_q}}\norm{\sum_{\alpha\in \Lambda(m,n)}\norm{v_{\alpha}}z^{\alpha}}=\norm{\sum_{\alpha\in \Lambda(m,n)}\norm{v_{\alpha}}z^{\alpha}}_{\mathcal{P}(^m\ell^n_q,X)}\\ &
	\leq \chi_{\text{mon}}(\mathcal{P}(^m\ell^n_q,X))\norm{\sum_{\alpha\in \Lambda(m,n)}v_{\alpha}z^{\alpha}}_{\mathcal{P}(^m\ell^n_q,X)}.
\end{align*} 
Therefore, 
\begin{equation*}
	\frac{1}{(K^n_m(B_{\ell^n_q},X,\lambda))^m}\leq \frac{1}{\lambda}\cdot \chi_{\text{mon}}(\mathcal{P}(^m\ell^n_q,X)),
\end{equation*}
and consequently shows that
\begin{equation}\label{P7-e-006}
	K^n_m(B_{\ell^n_q},X,\lambda)\geq \frac{\sqrt[m]{\lambda}}{\sqrt[m]{	\chi_{\text{mon}}(\mathcal{P}(^m\ell^n_q,X))}}.
\end{equation}
By the virtue of \eqref{P7-e-005} and \eqref{P7-e-006}, we obtain our desired equality \eqref{P7-e-015}. This completes the proof.
\end{pf}

The following result demonstrates a crucial link between the Bohr radius and the unconditional basis constant of spaces of vector-valued homogeneous polynomials.

\begin{thm}\label{P7-thm-02}
	Let $X$ be a finite dimensional complex Banach space and the canonical basis $\{e_{j}:1\leq j\leq n\}$ has unconditional basis constant $\chi(\{e_j\}^n_{j=1})=1$. Then for any $\lambda\geq 1$ and $1\leq q\leq \infty$ we have
	\begin{equation*}
		\frac{1}{3}\,\frac{\sqrt[m]{\lambda}}{\underset{m\in \mathbb{N}}{\sup}\sqrt[m]{	\chi_{\text{mon}}(\mathcal{P}(^m\ell^n_q,X))}} \leq K^n(B_{\ell^n_q},X,\lambda)\leq \min\left\{K(\mathbb{D},X,\lambda), \frac{\sqrt[m]{\lambda}}{\underset{m\in \mathbb{N}}{\sup}\sqrt[m]{	\chi_{\text{mon}}(\mathcal{P}(^m\ell^n_q,X))}}\right\}.
	\end{equation*}
\end{thm}

\begin{pf}
	The proof of the second inequality is just a consequence of the facts that 
	$K^n(B_{\ell^n_q},X,\lambda)\leq K(\mathbb{D},X,\lambda)$ for all $n$ and all $\lambda \geq 1$, and for all $m$
	\begin{equation*}
		K^n(B_{\ell^n_q},X,\lambda)\leq K^n_m(B_{\ell^n_q},X,\lambda)= \frac{\sqrt[m]{\lambda}}{\sqrt[m]{	\chi_{\text{mon}}(\mathcal{P}(^m\ell^n_q,X))}}.
	\end{equation*}
	Conversely, for the first inequality, let us assume that $f(z)=\sum_{\alpha \in \mathbb{N}^{n}_{0}}x_{\alpha}(f)z^{\alpha}$ be a bounded holomorphic function on $B_{\ell^n_q}$ with values in $X$. We want to prove first that for all $m$ and $z\in B_{\ell^n_q}$,
	\begin{equation}
		\norm{\sum_{\alpha\in \Lambda(m,n)}x_{\alpha}(f)z^{\alpha}}_{\mathcal{P}(^m\ell^n_q,X)}\leq 2(\norm{f}-\norm{x_{0}(f)}).
	\end{equation}
	 Given such $m$ and $z\in B_{\ell^n_q}$, we define for $\zeta\in \mathbb{D}$,
	 \begin{equation*}
	 	g(\zeta):= f(\zeta z)=\sum_{\alpha \in \mathbb{N}^{n}_{0}}x_{\alpha}(f)z^{\alpha}\zeta^{|\alpha|}=x_{0}
	 	(f)+\sum_{m=1}^{\infty}\left(\sum_{\alpha\in \Lambda(m,n)}x_{\alpha}(f)z^{\alpha}\right)\zeta^{m}.
	 \end{equation*}
 Then, $g$ is a vector-valued holomorphic function on $\mathbb{D}$ with $m$-th coefficient $\sum_{\alpha\in \Lambda(m,n)}x_{\alpha}(f)z^{\alpha}$ and $\norm{g(\zeta)}\leq 1$ for all $\zeta\in \mathbb{D}$. For any $\phi\in X^*$ with $\norm{\phi}\leq 1$, let us consider the function $\phi \circ g: \mathbb{D}\rightarrow \overline{\mathbb{D}}$ defined by
 \begin{equation*}
 (	\phi\circ g)(u)=\phi(x_0(f))+\sum_{m=1}^{\infty}\left(\sum_{\alpha\in \Lambda(m,n)}\phi(x_{\alpha}(f))z^{\alpha}\right)u^{m}.
 \end{equation*}
 Clearly, $\phi \circ g$ is holomorphic function on open unit disk $\mathbb{D}$ in $\mathbb{C}$ with the $m$-th coefficient
 \begin{equation*}
 	\phi\left(\sum_{\alpha\in \Lambda(m,n)}x_{\alpha}(f)z^{\alpha}\right)=\sum_{\alpha\in \Lambda(m,n)}\phi(x_{\alpha}(f))z^{\alpha}
 \end{equation*}
  and $|\phi(g(u))|\leq 1$ on $\mathbb{D}$. Therefore, $\text{Re}\,(\norm{f}-e^{i\theta}\phi(g))\geq 0$, where $\theta \in \mathbb{R}$ is such that $e^{i\theta}\phi(g(0))=|\phi(g(0))|$ for all $\phi\in B_{X^*}$. Using the Carath\'{e}odory's inequality \cite[Lemma 3.2]{aizenberg-2001} for all $m\in \mathbb{N}$, we have that
 \begin{equation*}
 	\bigg|\sum_{\alpha\in \Lambda(m,n)}\phi(x_{\alpha}(f))z^{\alpha}\bigg|\leq 2\,\text{Re} \, (\norm{f}-e^{i\theta}\phi(g(0)))=2(\norm{f}-|\phi(x_0(f))|)
 \end{equation*}
holds for any $\phi\in B_{X^*}$. Therefore, it follows that 
\begin{equation*}
	\norm{\sum_{\alpha\in \Lambda(m,n)}x_{\alpha}(f)z^{\alpha}}\leq 2\left(\norm{f}-\norm{x_0(f)}\right).
\end{equation*}
Hence, in view of Lemma \ref{P7-lem-03}, we obtain that for all $m\geq 1$ and all $z\in B_{\ell^n_q}$ 
\begin{align*}
		\sum_{\alpha\in \Lambda(m,n)}\norm{x_{\alpha}(f)\left(\frac{\lambda^{1/m}z}{3 \sqrt[m]{\chi_{\text{mon}}(\mathcal{P}(^m\ell^n_q,X))}}\right)}&\leq \frac{\lambda}{3^m}\norm{\sum_{\alpha\in \Lambda(m,n)}x_{\alpha}(f)z^{\alpha}}\\&
		\leq \frac{\lambda}{3^m}\cdot 2\left(\norm{f}-\norm{x_0(f)}\right).
\end{align*}
So, for each
\begin{equation*}
	z\in \left(\frac{\lambda^{1/m}}{3 \sqrt[m]{\chi_{\text{mon}}(\mathcal{P}(^m\ell^n_q,X))}}\right)B_{\ell^n_q}
\end{equation*}
we have
\begin{equation*}
	\sum_{\alpha \in \mathbb{N}^{n}_{0}}\norm{x_{\alpha}(f)z^{\alpha}}=\norm{x_0(f)}+\sum_{m=1}^{\infty}\norm{x_{\alpha}(f)z^{\alpha}}\leq \norm{x_0(f)}+2\left(\norm{f}-\norm{x_0(f)}\right)\lambda\sum_{m=1}^{\infty}\frac{1}{3^m}=\lambda\norm{f},
\end{equation*}
which shows that
\begin{equation*}
\frac{1}{3}\,\frac{\lambda^{1/m}}{ \sqrt[m]{\chi_{\text{mon}}(\mathcal{P}(^m\ell^n_q,X))}}\leq K^n(B_{\ell^n_q},X,\lambda).
\end{equation*}
This completes the proof.
 \end{pf}
\section{Estimate for Bohr radius}\label{P7-sec-04}
The following lemma is required in order to prove our main result. It is a consequence of a classical result by Wiener \cite[Lemma 8.4, p. 183]{defant-book}.
\begin{lem}\label{P7-lem-01}
	Let $f(z)=\sum_{\alpha \in \mathbb{N}^{n}_{0}}x_{\alpha}z^{\alpha}$ be a bounded holomorphic functions with values in a complex Banach space  $X$ such that $\sup_{z\in B_{\ell^n_q}}\norm{f(z)}_{X}\leq 1.$ Then the following holds 
\begin{equation*}
	\sup_{z\in \partial B_{\ell^n_q}} \sum_{\alpha \in \Lambda(m,n)}\norm{x_{\alpha}z^{\alpha}}\leq \widetilde{S}(m,n) (1-\norm{x_0}^2)
\end{equation*}
for $m\geq 1$.
\end{lem}
\begin{pf}
	Let $z\in \partial B_{\ell^n_q}$ be fixed. For $u\in \mathbb{D}$, we consider the function 
	\begin{equation*}
		g(u):=f(uz)=x_{0}+\sum_{m=1}^{\infty} \left(\sum_{\alpha\in \Lambda(m,n)}x_{\alpha}z^{\alpha}\right)u^m.
	\end{equation*}
Then $g$ is a holomorphic function on unit disk $\mathbb{D}$ into the complex Banach space $X$ with the $m$-th coefficient $\sum_{\alpha\in \Lambda(m,n)}x_{\alpha}z^{\alpha}$ and $\norm{g(u)}\leq 1$ for all $u\in \mathbb{D}$. Consequently, we note that for any $\phi\in X^*$ with $\norm{\phi}\leq 1$, the function $\phi \circ g: \mathbb{D} \rightarrow \overline{\mathbb{D}}$ is holomorphic with the $m$-th coefficient 
\begin{equation*}
	\phi\left(\sum_{\alpha\in \Lambda(m,n)}x_{\alpha}z^{\alpha}\right)=\sum_{\alpha\in \Lambda(m,n)}\phi(x_{\alpha})z^{\alpha}.
\end{equation*}
Therefore, by the virtue of the Wiener's result \cite[Lemma 8.4, p. 183]{defant-book}, we obtain
\begin{equation*}
\bigg|\phi\left(\sum_{\alpha\in \Lambda(m,n)}x_{\alpha} z^{\alpha}\right)\bigg|\leq 1-|\phi(x_0)|^2
\end{equation*}
holds for any $\phi \in B_{X^*}$. Hence, it shows that 
\begin{equation*}
	\sup_{z\in \partial B_{\ell^n_q}} \norm{ \sum_{\alpha\in \Lambda(m,n)}x_{\alpha}z^{\alpha}}\leq 1-\norm{x_{0}}^2.
\end{equation*}
Using the Definition \ref{defn-sidon}, it follows that 
\begin{equation*}
	\sup_{z\in \partial B_{\ell^n_q}}\sum_{\alpha\in \Lambda(m,n)}\norm{x_{\alpha}z^{\alpha}} \leq \widetilde{S}(m,n)\sup_{z\in \partial B_{\ell^n_q}}\norm{ \sum_{\alpha\in \Lambda(m,n)}x_{\alpha}z^{\alpha}}\leq \widetilde{S}(m,n) (1-\norm{x_{0}}^2).
\end{equation*}
This completes the proof.
\end{pf}

By the virtue of Corollary \ref{P7-cor-01}, we obtain for $n\geq 2$ the function
\begin{equation}\label{P7-e-014}
	\widetilde{H_n}(x):=\sum_{m=1}^{\infty}\widetilde{S}(m,n)x^m
\end{equation}
is well-defined and differentiable in $(-1,1)$, and so in $[0,1)$. Analogously to Lemma \ref{P7-lem-04}, we can say that $\widetilde{H_n}$ is strictly increasing on $[0,1)$ and $\widetilde{H_n}(x)>x$ for all $x\in (0, 1)$, whereas $\widetilde{H_n}(0)=0$. By intermediate mean value theorem, for every $M>0$, there exists a unique $R\in (0,1)$ such that $\widetilde{H_n}(R)=M$. Hence, for any $1\leq \lambda<\infty$, by considering $M=\lambda/2$, the following definition makes sense. 

\begin{defn}
	For every $n\in \mathbb{N}$ and $1\leq \lambda<\infty$, we denote by $\gamma_n\in (0,1)$ the unique solution of the following equation:
	\begin{equation}\label{P7-e-001-gamma}
		x + \sum_{m=2}^{\infty}\widetilde{S}(m,n) x^m=\frac{\lambda}{2}.
	\end{equation}
\end{defn}

Next, we prove our main theorem which shows that $\gamma_{n}$ is a lower bound of Bohr radius $K^n(B_{\ell^n_q},X,\lambda)$ and it is better that $\beta_n$. Further, we prove that our lower estimate is the improved bound than that of from \cite{defant-2012}.
\begin{thm}\label{P7-thm-01}
	Let $X$ be any finite dimensional complex Banach space and $n\in \mathbb{N}$ with $n\geq 2$. For $\lambda>1$ and $1\leq q< \infty$ suppose $K^n(B_{\ell^n_q},X,\lambda)$ is the $n$-dimensional Bohr radius and the numbers $\beta_n$, $\gamma_n$ that had been defined by \eqref{P7-e-001-beta} and \eqref{P7-e-001-gamma}. Then we have
	\begin{equation*}
		\beta_n<\gamma_n\leq K^n(B_{\ell^n_q},X,\lambda)\leq \frac{\gamma_n}{K(\mathbb{D},X,\lambda)}.
	\end{equation*}
\end{thm}

\begin{pf}
	In view of Corollary \ref{P7-cor-01}, we obtain $\widetilde{H_n}(x)\leq H_n(x)$ for all $x\in (0,1)$, where $H_n$ and $\widetilde{H_n}$ are defined by \eqref{P7-e-012} and \eqref{P7-e-014}, respectively. We observe that 
	\begin{equation*}
		H_n(\beta_n)=\frac{\lambda}{2}=\widetilde{H_n}(\gamma_n)<H_n(\gamma_n). 
	\end{equation*}
Since $H_n$ is strictly monotonic (due to Lemma \ref{P7-lem-04}) for all $n\in \mathbb{N}$, it follows that $\beta_n<\gamma_n$.\\

Next we want to prove that $\gamma_n\leq K^n(B_{\ell^n_q},X,\lambda)$. Without loss of generality, let us assume that $f(z)=\sum_{\alpha \in \mathbb{N}^{n}_{0}}x_{\alpha}z^{\alpha}$ is a bounded holomorphic function defined on unit ball of $\ell^n_q$ space with 
\begin{equation*}
	\sup_{z\in B_{\ell^n_q}}\norm{f(z)}_{X}\leq 1.
\end{equation*}
We choose a point $z=(z_1,\ldots, z_n)\in B_{\ell^n_q}$ with $|z_k|\leq \gamma_n$ for all $k \in \{1,2,\ldots, n\}$.
By the virtue of Lemma \ref{P7-lem-01}, we obtain
\begin{align*}
	\sup_{z\in \gamma_n B_{\ell^n_q}}\sum_{\alpha \in \mathbb{N}^{n}_{0}}\norm{x_{\alpha}z^{\alpha}}&\leq \norm{x_0} + \sum_{m=1}^{\infty}\left(\sup_{z\in B_{\ell^n_q}}\sum_{\alpha\in \Lambda(m,n)} \norm{x_{\alpha}}|z^{\alpha}|\right)\left(\gamma_{n}\right)^m\\ &
	\leq \norm{x_0} + \sum_{m=1}^{\infty} (1-\norm{x_0}^2) \widetilde{S}(m,n)\gamma^m_{n}\\ &
	= \norm{x_0}+(1-\norm{x_0}^2) \left(\sum_{m=1}^{\infty}\widetilde{S}(m,n)\gamma^m_{n}\right)\\ &
	= \norm{x_0} + (1-\norm{x_0}^2)\cdot \frac{\lambda}{2}\\ &
	\leq \norm{x_{0}} + 2(1-\norm{x_0})\cdot \frac{\lambda}{2}=\lambda.
\end{align*}
Therefore, it follows from the definition of $K^n(B_{\ell^n_q},X,\lambda)$ that $\gamma_{n}\leq K^n(B_{\ell^n_q},X,\lambda)$.\\

Next our aim is to show that $K^n(B_{\ell^n_q},X,\lambda)\leq \gamma_{n}/K(\mathbb{D},X,\lambda)$. To prove this, we first fix $n\in \mathbb{N}$ with $n\geq 2$. On studying the related research on Bohr radius for vector-valued holomorphic functions, the authors in \cite{Subhadip-Vasu-P3} have proved that $K^n(B_{\ell^n_q},X,\lambda)\leq \left(K(\mathbb{D},X,\lambda)\right)^{1/q}$ for every $q\in [1,\infty)$. 
 It is worth noting that if $K^n(B_{\ell^n_q},X,\lambda)=(K(\mathbb{D},X,\lambda))^{1/q}$, then there is noting to prove because we obtain the exact the value of $K^n(B_{\ell^n_q},X)$ in this case. Hence, we consider $K^n(B_{\ell^n_q},X,\lambda)<(K(\mathbb{D},X,\lambda))^{1/q}$ for $1\leq q<\infty$ to avoid such trivial situation. Let us fix
 \begin{equation*}
 	\delta \in \left(\frac{1}{K(\mathbb{D},X,\lambda)^{1/q}}, \frac{1}{K^n(B_{\ell^n_q},X,\lambda)}\right).
 \end{equation*}
  Then we have $\delta>1/K(\mathbb{D},X,\lambda)^{1/q}$. We assume $\zeta=\delta^q$ for our convenience. An easy computation shows that 
  \begin{equation*}
  	\frac{1}{\zeta}<1-\frac{\frac{1}{K(\mathbb{D},X,\lambda)}-1}{\zeta},
  \end{equation*}
  and thus we can choose 
  \begin{equation}\label{P7-e-002}
  	r\in \left(\frac{1}{\zeta}, 1-\frac{\frac{1}{K(\mathbb{D},X,\lambda)}-1}{\zeta}\right).
  \end{equation}

   For each $m\in \mathbb{N}$, the space $\mathcal{P}(^m\ell^n_q,X)$ is finite dimensional, so we can find $Q_m\in \mathcal{P}(^m\ell^n_q,X)$ satisfying $\norm{Q_m}_{\infty}=1$ and $\norm{Q_m}_{1}=\widetilde{S}(m,n)$.
  For each $\phi\in X^*$ with $\norm{\phi}\leq 1$, the series 
  \begin{equation}\label{P7-e-001}
  	\sum_{m=1}^{\infty}\frac{r^m}{2} Q_m(z),\quad r\in (0,1),
  \end{equation}
converges absolutely and uniformly on $B_{\ell^n_q}$. Therefore, the series expansion \eqref{P7-e-001} defines a holomorphic function $D(z)$ for $z\in B_{\ell^n_q}$. Further, we observe that $	\sup_{\phi \in B_{X^*}} \norm{\phi(Q_m(z))}_{\infty}=1.$ Indeed, we obtain the following:
\begin{equation*}
\sup_{\phi \in B_{X^*}}	|\phi(D(z))|\leq\sup_{\phi \in B_{X^*}} \sum_{m=1}^{\infty}\frac{r^m}{2} \norm{\phi(Q_m(z))}_{\infty} =\frac{r}{2(1-r)} \quad \mbox{for all}\,\,\, z\in B_{\ell^n_q}
\end{equation*}
gives 
\begin{equation*}
	\norm{D}_{\infty}\leq \frac{r}{2(1-r)}.
\end{equation*}
Let $D(z):=\sum_{\alpha \in \mathbb{N}^{n}_{0}}x_{\alpha}z^{\alpha}$ be the monomial series expansion of $D$ on $B_{\ell^n_q}$. Let us denote $R_n:=\zeta \gamma_{n}$. It is easy to observe that $0<R_n<\zeta K^n(B_{\ell^n_q},X,\lambda)$. In view of \eqref{P7-e-002}, we have $\zeta r>1$ and 
\begin{equation*}
	\zeta >\frac{\frac{1}{K(\mathbb{D},X,\lambda)}-1}{1-r}.
\end{equation*}
Therefore, for any $z\in \zeta \gamma \partial B_{\ell^n_q}$ we obtain the following:
\begin{align*}
	\sum_{\alpha \in \mathbb{N}^{n}_{0}}\norm{x_{\alpha}z^{\alpha}}&=\sum_{m=1}^{\infty}\sum_{\alpha\in \Lambda(m,n)} \norm{x_{\alpha}z^{\alpha}}=\sum_{m=1}^{\infty}\sum_{\alpha\in \Lambda(m,n)} \norm{x_{\alpha}z^{\alpha}} (\zeta \gamma_{n})^m\\ &
	=\sum_{m=1}^{\infty}(\zeta \gamma_{n})^m \norm{Q_m}_{1} r^m=\sum_{m=1}^{\infty}\zeta^m \gamma^m_{n} r^m \widetilde{S}(m,n)\\&
	=\zeta r\sum_{m=1}^{\infty}\widetilde{S}(m,n)(\gamma_{n})^m(\zeta r)^{m-1}\\&
	> \zeta r \sum_{m=1}^{\infty}\widetilde{S}(m,n)(\gamma_{n})^m\\ &
	=\zeta r \cdot \frac{\lambda}{2}
	>\frac{\frac{1}{K(\mathbb{D},X,\lambda)}-1}{1-r}\cdot \frac{r\lambda}{2}>\frac{r\lambda}{2(1-r)}\geq \lambda\norm{D}_{\infty}.
\end{align*}
In view of the definition of $K^n(B_{\ell^n_q},X,\lambda)$, it follows that $\zeta \gamma_{n}=R_{n}\geq K^n(B_{\ell^n_q},X,\lambda)$. Now, letting $\zeta \rightarrow (1/K(\mathbb{D},X,\lambda))^{+}$, we obtain that $K^n(B_{\ell^n_q},X,\lambda)\leq \gamma_{n}/K(\mathbb{D},X,\lambda)$. This completes the proof.
\end{pf}

Now, our final goal is to show that for every $n\geq 2$, $\gamma_{n}$ is an improved lower bound than that from \cite{defant-2012}. From the definition of unconditional basis constant of homogeneous polynomials $\chi_{\text{mon}}(\mathcal{P}(^m\ell^n_q,X))$, it is evident that 
\begin{equation*}
	\chi_{\text{mon}}(\mathcal{P}(^m\ell^n_q,X))=\widetilde{S}(m,n).
\end{equation*}
Let us define for each $n\in \mathbb{N}$ 
\begin{equation}\label{P7-e-007}
	\Gamma_{n}:=\sup\left\{\widetilde{S}(m,n)^{1/m}:m\in \mathbb{N}\right\}.
\end{equation}
As a consequence of Theorem \ref{P7-thm-02}, we have $\sqrt[m]{\lambda}/(3\Gamma_{n})\leq K^n(B_{\ell^n_q},X,\lambda)$, $\lambda\geq 1$. 

\begin{prop}
Let $X$ be a finite dimensional complex Banach space and $\lambda>1$. Then, for any $1\leq q\leq \infty$ we have
\begin{equation*}
	C\left(\frac{\log n}{n}\right)^{1-\frac{1}{\min\{q,2\}}}<\gamma_{n}\leq K^n(B_{\ell^n_q},X,\lambda),
\end{equation*}
where $C$ being the constant depending only on $\lambda$ and $X$.
\end{prop}

\begin{pf}
	The right hand inequality has already been proved in Theorem \ref{P7-thm-01}. To prove the left hand inequality, let us fix $n\geq 2$. Firstly, our aim is to show that there exists $m\in \mathbb{N}$ such that 
	\begin{equation*}
		\widetilde{S}(m,n)^{1/m}<\Gamma_{n}, \quad n\geq 2,
	\end{equation*}
where $\Gamma_{n}$ is given by \eqref{P7-e-007}. Indeed, if this is not true, we suppose that $\widetilde{S}(m,n)=\Gamma^m_{n}$ for all $m\in \mathbb{N}$. For the case $m=1$, we want to show that $\widetilde{S}(1,n)=1$ for all $n\in \mathbb{N}$. Indeed, from the definition of $\widetilde{S}(m,n)$, we have already that $\widetilde{S}(1,n)\geq 1$. Let $Q(z)=\sum_{j=1}^{n}x_j z_j$ be a $1$-homogeneous polynomial with $x_j\in X$ and $z=(z_1,\ldots,z_n)\in \ell^n_q$. For every $\phi\in B_{X^*}$, consider the function 
\begin{equation*}
	(\phi\circ Q)(z)=\phi\left(\sum_{j=1}^{n}x_j z_j\right)=\sum_{j=1}^{n}\phi(x_j)z_j.
\end{equation*}
Let $\theta_{j}\in [0,2\pi]$ be such that $|\phi(x_j)z_j|=\phi(x_j)z_je^{i\theta_{j}}$ for $j=1,2,\ldots,n$. Now, we have that
\begin{align*}
	\sum_{j=1}^{n}\norm{x_j z_j}&=\sup_{\phi \in B_{X^*}}\sum_{j=1}^{n}|\phi(x_j)z_j|=\sup_{\phi \in B_{X^*}}\sum_{j=1}^{n}\phi(x_j)z_je^{i\theta_{j}}\\ &
	=\sup_{\phi \in B_{X^*}} \phi\left(\sum_{j=1}^{n}x_jz_je^{i\theta_{j}}\right)=\sup_{\phi \in B_{X^*}}\phi\left(Q(z_1e^{i\theta_{1}},\ldots,z_ne^{i\theta_{n}})\right)\leq \norm{Q}_{\infty}.
\end{align*}
Therefore, it follows that $\widetilde{S}(1,n)\leq 1$. Consequently, $\Gamma_{n}=1$, and so $\widetilde{S}(m,n)=1$ for all $m\in \mathbb{N}$, which is a contradiction. Hence, we have proved our claim.\\

\noindent
The function $\widetilde{H_n}(x)=\sum_{m=1}^{\infty}\widetilde{S}(m,n)x^m$ satisfies that $\widetilde{H_n}(\gamma_{n})=\lambda/2$, and 
\begin{equation*}
	\widetilde{H_n}\left(\frac{1}{3\Gamma_{n}}\right)=\sum_{m=1}^{\infty}\widetilde{S}(m,n)\cdot \frac{1}{3^m \Gamma^m_{n}}=\sum_{m=1}^{\infty}\frac{\widetilde{S}(m,n)}{\Gamma^m_{n}}\cdot \frac{1}{3^m}<\sum_{m=1}^{\infty}\frac{1}{3^m}=\frac{1}{2}<\frac{\lambda}{2}=\widetilde{H_n}(\gamma_{n}).
\end{equation*}
Since $\widetilde{H_n}$ is strictly increasing, we obtain $1/(3\Gamma_{n})<\gamma_{n}$. Finally, in view of the result \cite[Theorem 1.1]{defant-2011-angew}, we have 
\begin{equation*}
		C\left(\frac{\log n}{n}\right)^{1-\frac{1}{\min\{q,2\}}}\leq \frac{1}{3\Gamma_{n}},
\end{equation*}
where $C$ depends on $X$ and $\lambda$. This completes the proof.
\end{pf}

\vspace{3mm}

\noindent\textbf{Compliance of Ethical Standards:}\\

\noindent\textbf{Conflict of interest.} The authors declare that there is no conflict  of interest regarding the publication of this paper.
\vspace{1.5mm}

\noindent\textbf{Data availability statement.}  Data sharing is not applicable to this article as no datasets were generated or analyzed during the current study.
\vspace{1.5mm}

\noindent\textbf{Authors contributions.} Both the authors have made equal contributions in reading, writing, and preparing the manuscript.\\

\noindent\textbf{Acknowledgment:} 
%The authors would like to express their sincerest gratitude to the referee for careful reading of the manuscript and many valuable suggestions, which greatly helped to improve the clarity of the exposition in this manuscript. 
%The research of the first named author is supported by SERB-CRG (DST), Govt. of India.  
The research of the second named author is supported by DST-INSPIRE Fellowship (IF 190721),  New Delhi, India.


\begin{thebibliography}{99}
	

	
	
	\bibitem{aizn-2000a} {\sc L. Aizenberg}, Multidimensional analogues of Bohr's theorem on power series, \textit{Proc. Amer. Math. Soc.} {\bf 128} (2000), 1147--1155.
	
	\bibitem{aizn-2000b} {\sc L. Aizenberg, A. Aytuna}, and {\sc P. Djakov}, An abstract approach to Bohr's phenomenon, {\it Proc. Amer. Math. Soc.} {\bf 128} (2000), 2611--2619.
	
	
	\bibitem{aizenberg-2001} {\sc L. Aizenberg, A. Aytuna}  and {\sc P. Djakov}, Generalization of theorem on Bohr for bases in spaces of holomorphic functions of several complex variables, 
	{\it J. Math. Anal. Appl.} {\bf  258} (2001), 429--447.
	
	\bibitem{aizn-2005} {\sc L. Aizenberg}, Generalization of Carath\'{e}odory’s inequality and the Bohr radius for multidimensional power series, Selected Topics in Complex Analysis. \textit{Operator Theory: Advances and Applications}, {\bf 158} (2005), Birkhäuser Basel, 87--94.
	
	
%	\bibitem{aizn-2007} {\sc L. Aizenberg}, Generalization of results about the Bohr radius for power series, {\it Stud. Math.}  {\bf 180}  (2007), 161--168.  
	
%	\bibitem{aizenberg-2012} {\sc L. Aizenberg}, Remarks on the Bohr and Rogosinski phenomena for power series, {\it Anal. Math. Phys.} {\bf 2} (2012), 69--78.
	

	

	
%	\bibitem{Subhadip-Vasu-P1} {\sc V. Allu, H. Halder,} and {\sc S. Pal}, Bohr and Rogosinski inequalities for operator valued holomorphic functions, {\it Bull. Sci. Math.} {\bf 182} (2023), 103214, 18 pp.
	
	
	\bibitem{Subhadip-Vasu-P3} {\sc V. Allu, H. Halder,} and {\sc S. Pal}, Multidimensional Bohr radii for holomorphic functions with values in complex Banach spaces, see https://arxiv.org/abs/2308.07825 (2023).
	
%	\bibitem{Himadri-Vasu-canad-2023} {\sc V. Allu } and {\sc H. Halder}, Bohr operator on operator-valued polyanalytic functions on simply connected domains {\it Canad. Math. Bull.}, (2023), 12 pp, doi:10.4153/S0008439523000541.
	
	\bibitem{bala-studia-2006} {\sc R. Balasubramanian, B. Calado}, and {\sc H. Queff\'{e}lec}, The Bohr inequality for ordinary Dirichlet series, {\it Studia Math.} {\bf 175} (2006), 285--304. 
	
%	\bibitem{anderson-2006} {\sc J. M. Anderson} and {\sc J. Rovnyak}, On Generalized Schwarz-Pick Estimates, {\it Mathematika} {\bf 53} (2006), 161--168.
	
	
	%	\bibitem{bene-2004} {\sc C. B{\rm $\acute{E}$}n$\acute{E}$teau, A. Dahlner} and {\sc D. Khavinson}, Remarks on the Bohr phenomenon, {\it  Comput. Methods Funct. Theory}  {\bf 4}  (2004), 1--19.
	
	%   \bibitem{bharandar-2014} {\sc SV Bharanedhar} and {\sc S. Ponnusamy}, Coefficient conditions for harmonic univelent mappings and hypergeometric mappings, {\it Rocky Mountain J. Math.}	{\bf 44} (2014) 753--777.
	
%	\bibitem{Ayt & Dja & BLMS & 2013} {\sc A. Aytuna} and {\sc P. Djakov}, Bohr property of bases in the space of entire functions and its generalizations, Bull. London Math. Soc. \textbf{45}(2)(2013), 411--420.
	
	\bibitem{bayart-advance-2014} {\sc F. Bayart, D. Pellegrino}, and {\sc J. B. Seoane-Sep$\rm \acute{U}$lveda}, The Bohr radius of the $n$-dimensional polydisk is equivalent to $\sqrt{(\log n)/n}$, {\it Adv. Math.} {\bf 264} (2014), 726--746.
	
	%	\bibitem{beardon-minda-hyperbolic-density} {\sc A. F. Beardon} and {\sc D. Minda}, The hyperbolic metric and geometric function theory, {\it Proceedings of the International Workshop on Quasiconformal Mappings and their Applications (IWQCMA05)}.
	
	\bibitem{Bernal-2021} {\sc L. Bernal-Gonz\'{a}lez, H. J. Cabana, D. Garcia, M. Maestre, G. A. Mu\~{n}oz-Fern\'{a}ndez} and {\sc J. B. Seoane-Sep\'{u}iveda}, A new approach towards estimating the $n$-dimensional Bohr radius, {\it RACSAM} {\bf 115} (2021), 44.
	
%	\bibitem{bene-2004} {\sc C. B{\'e}n{\'e}teau}, {\sc A. Dahlner} and {\sc D. Khavinson}, Remarks on the Bohr phenomenon, {\it Comput. Methods Funct. Theory} \textbf{4}(1) (2004), 1-19.
	

	
	\bibitem{Blasco-OTAA-2010} {\sc O. Blasco}, The Bohr radius of a Banach space, {\it In Vector measures, integration and related topics, 5964, Oper. Theory Adv. Appl., 201, Birkh{\"a}user Verlag, Basel, 2010.}
	
	\bibitem{Blasco-Collect-2017} {\sc O. Blasco}, The $p$-Bohr radius of a Banach space, {\it Collect. Math.} {\bf 68} (2017), 87--100.
	
	
	\bibitem{boas-1997} {\sc H.P. Boas} and {\sc D. Khavinson}, Bohr's power series theorem in several variables, {\it Proc. Amer. Math. Soc.}  {\bf 125} (1997), 2975--2979.
	
		\bibitem{boas-2000} {\sc H.P. Boas}, Majorant Series,  {\it J. Korean Math. Soc.}  {\bf 37} (2000), 321--337.
	
	\bibitem{Bohr-1914} {\sc H. Bohr}, A theorem concerning power series,  {\it Proc. Lond. Math. Soc.} s2-13 (1914), 1--5.
	
%	\bibitem{bombieri-1962} {\sc E. Bombieri}, Sopra un teorema di H. Bohr e G. Ricci sulle funzioni maggioranti delle serie di potenze, {\it Boll. Un. Mat. Ital.} {\bf 17} (1962), 276--282.
	
%	\bibitem{bombieri-2004} {\sc E. Bombieri} and {\sc J. Bourgain}, A remark on Bohr's inequality, {\it Internat. Math. Res. Notices} {\bf 80} (2004), 4307--4330.
	

	
%	\bibitem{das-2023} {\sc N. Das}, Estimates for generalized Bohr radii in one and higher dimensions, {\it Canad. Math. Bull.} {\bf 66} (2023), 682--699.
	
	\bibitem{defant-2003} {\sc A. Defant, D. Garc\'{i}a}, and {\sc  M. Maestre}, Bohr power series theorem and local Banach space theory, {\it J. Reine Angew. Math.} {\bf 557} (2003), 173--197.
	
	\bibitem{defant-book} {\sc A. Defant, D. García,  M. Maestre}, and {\sc P. Sevilla-Peris}, \textit{Dirichlet Series and Holomorphic Functions in High Dimensions}, New Mathematical Monographs: {\bf 37}, Cambridge University Press, Cambridge, (2019).
	
%	\bibitem{defant-2004} {\sc A. Defant, D. Garc\'{i}a}, and {\sc  M. Maestre}, Estimates for the first and second Bohr radii of Reinhardt domains, {\it J. Appr. Theory} {\bf 128} (2004), 53--68.
	
	\bibitem{defant-2006} {\sc A. Defant} and {\sc L. Frerick}, A logarithmic lower bound for multi-dimenional Bohr radii, {\it Israel J. Math.} {\bf 152} (2006), 17--28.
	
%	\bibitem{defant-2007} {\sc A. Defant,  M. Maestre}, and {\sc  C. Prengel}, The arithmetic Bohr radius, {\it Q. J. Math.} {\bf 59} (2008), 189--205.
	
	
	\bibitem{defant-2011} {\sc A. Defant, L. Frerick, J. Ortega-Cerd${\rm \grave{A}}$, M. Ouna{\"i}es}, and {\sc K. Seip}, The Bohnenblust-Hille inequality for homogeneous polynomils in hypercontractive, {\it Ann. of Math.} {\bf 174} (2011), 512--517.

	\bibitem{defant-2011-angew} {\sc A. Defant} and {\sc L. Frerick}, The Bohr radius of the unit ball of $\ell^n_p$, {\it J. Reine Angew Math.} {\bf 660} (2011), 131--147.

	\bibitem{defant-2012} {\sc A. Defant, M. Maestre}, and {\sc U. Schwarting}, Bohr radii for vector valued holomorphic functions, {\it Adv. Math.} {\bf 231} (2012), 2837--2857.
	
%	\bibitem{defant-book} {\sc A. Defant, D. García,  M. Maestre}, and {\sc P. Sevilla-Peris}, \textit{Dirichlet Series and Holomorphic Functions in High Dimensions} (New Mathematical Monographs) Cambridge: Cambridge University Press, (2019) doi:10.1017/9781108691611.

\bibitem{defant-2018} {\sc A. Defant, M. Masty\l{o}}, and {\sc S. Schl\"{u}ters}, On Bohr radii of finite dimensional complex Banach spaces, {\it Funct. Approx. Comment. Math.} {\bf 59} (2018), 251--268.
	
	
%	\bibitem{Dixon & BLMS & 1995} {\sc P. G. Dixon}, Banach algebras satisfying the non-unital von Neumann inequality, \textit{Bull. Lond. Math. Soc.} \textbf{27} (4) (1995), 359--362.
	
%	\bibitem{Djakov & Ramanujan & J. Anal & 2000} {\sc P. B. Djakov} and {\sc M. S. Ramanujan}, A remark on Bohr's theorem and its generalizations, \textit{J. Anal.} \textbf{8} (2000), 65--77.
	
	\bibitem{Dineen-Timoney-1989} {\sc S. Dineen} and {\sc R. M. Timoney}, Absolute bases, tensor products and a theorem of Bohr, \textit{Studia Math.} \textbf{94} (1989), 227--234.
	
	
	\bibitem{kumar-2023-JMAA} {\sc S. Kumar} and {\sc R. Manna}, Revisit of multi-dimensional Bohr radius, {\it J. Math. Anal. Appl.} {\bf 523} (2023), 127023.
	
	\bibitem{kumar-2023-arxiv} {\sc S. Kumar} and {\sc R. Manna}, Multi-dimensional Bohr radii of Banach space valued holomorphic functions, see https://arxiv.org/pdf/2303.17416 (2023).

	\bibitem{paulsen-2002} {\sc Vern I. Paulsen, Gelu Popescu} and {\sc Dinesh Singh}, On Bohr's inequality, {\it Proc. Lond. Math. Soc.} s3-85 (2002), 493--512.
	
	%\bibitem{prengel-2005} {\sc C.Prengel}, {\it Domains of convergence in infinite dimensional holomorphy}, Ph.D. Thesis, University of Oldenburg, 2005.
	

	\bibitem{popescu-2019} {\sc G. Popescu}, Bohr inequalities for free holomorphic functions on polyballs, {\it Adv. Math.} {\bf 347} (2019), 1002-1053.
	
	

	%\bibitem{ronning-1993} {\sc F. Ronning}, Uniformly convex functions and a corresponding class of starlike functions, {\it Proc. Amer. Math. Soc.} {\bf 118} (1993), 189--196.
	
	%	\bibitem{sharma-2016} {\sc K. Sharma}, {\sc N. K. Jain} and {\sc V. Ravichandran}, Starlike functions associated with cardiod, {\it Afr. Mat.} {\bf 27} (2016), 923--939.
	
	%	\bibitem{sidon-1927} {\sc S. Sidon}, Uber einen satz von Hernn Bohr, {\it Math. Zeit.}  {\bf 26} (1927),  731-732.
	
	%\bibitem{silverman-1990} {\sc H. Silverman} and {\sc E.M. Silvia}, Subclasses of univalent functions starlike with respect to a boundary point, {\it Houston J. Math.} {\bf 16} (1990), 289--299. 
	
	%	\bibitem{rogosinski-1943} {\sc W. Rogosinski}, On the coefficients of subordinate functions, {\it Proc. London Math. Soc.} {\bf 48} (1943), 48--82.
	
	%	\bibitem{Ruscheweyh-1985} {\sc St. Ruscheweyh}, Two remarks on bounded analytic functions, \emph{Serdica} \textbf{11}(1) (1985), 731--732.
	
	
%	\bibitem{stempak-1994} {\sc K. Stempak}, Ces{\'a}ro averaging operators, {\it Proc. R. Soc. Edinb., Sect. A, Math.} {\bf 124} (1994), 121--126.
	
	%	\bibitem{Y. Sun} {\sc Y. Sun}, {\sc Y-P Jiyang} and {\sc A. Rasila}, On a certain subclass of close-to-convex harmonic mappings, {\it Complex Var. Elliptic Equ.} {\bf 61} (2016) 1627-1643.
	
	%	\bibitem{vasu-book} {\sc Derek K. Thomas}, {\sc Nikola Tuneski} and {\sc Allu Vasudevarao}, {\it Univalent functions}. A primer, De Gruyter Studies in Mathematics, {\bf 69}. De Gruyter, Berlin, 2018.
	
	%	\bibitem{tomic-1962} {\sc M. Tomic}, Sur un theoreme de H. Bohr,  {\it Math. Scand.}  {\bf 11} (1962), 103--106.
	
	%	\bibitem{wang-2006-a} {\sc Z. Wang}, {\sc C. Gao} and {\sc S. Yuan}, On certain subclass of close-to-convex functions, {\it Acta Math Acad Paedagog Nyhazi (N. S.)} {\bf 22} (2006), 171--177. (electronic)
	
	%\bibitem{Wang-2006-b} {\sc Z-G Wang}, {\sc C-Y Gao} and {\sc S-M Yuan}, On certain subclasses of close-to-convex and quasi-convex functions with respect to $k$-symmetric points, {\it J. Math. Anal. Appl.} {\bf 322} (2006), 97--106.
	
	%	\bibitem{Wani-2019} {\sc L. A. Wani} and {\sc A. Swaminathan}, Starlike and convex functions associated with nephroid domain, {\it Bull. Malays. Math. Sci. Soc.} doi:10.1007/s40840-020-00935-6. 2020
	
	%\bibitem{zhu-2016} {\sc J. Zhu}, Coefficients estimate for harmonic $\nu$-Bloch-mappings and harmonic $K$-quasiconformal mappings, {\it Bull. Malays. Math. Sci. Soc.} {\bf 39} (2016), 349--358.
\end{thebibliography}
\end{document}